\documentclass{amsart}

\DeclareMathSymbol{\twoheadrightarrow}  {\mathrel}{AMSa}{"10}

\def\Q{{\mathbb Q}}
\def\Z{{\mathbb Z}}
\def\C{{\mathbb C}}

\def\F{{\mathbb F}}
\def\P{{\mathbb P}}

     \def\CC{\mathfrak{C}}

\def\Gal{\mathrm{Gal}}
          \def\tr{\mathrm{tr}}
             \def\ord{\mathrm{ord}}

\def\End{\mathrm{End}}
\def\Aut{\mathrm{Aut}}

       \def\Lie{\mathrm{Lie}}

\def\II{\mathrm{Id}}

\def\cc{\mathfrak{c}}

\def\fchar{\mathrm{char}}
\def\GL{\mathrm{GL}}
\def\SL{\mathrm{SL}}

                      \def\ST{{\mathbf{S}}}

\def\dim{\mathrm{dim}}
                  \def\Fr{\mathrm{Fr}}

\def\O{{\mathcal O}}

       \def\p{\mathfrak{p}}
                        \def\sL{\mathfrak{sl}}
    \def\g{\mathfrak{g}}

\newtheorem{thm}{Theorem}[section]
\newtheorem{lem}[thm]{Lemma}

\newtheorem{prop}[thm]{Proposition}
\theoremstyle{definition}

\newtheorem{ex}[thm]{Example}
\newtheorem{rem}[thm]{Remark}

\title[Abelian varieties  without homotheties]
{Abelian varieties  without homotheties}
\author[Yuri\ G.\ Zarhin]{Yuri\ G.\ Zarhin}
\address{Department of Mathematics, Pennsylvania State University,
University Park, PA 16802, USA}
\address{Steklov Mathematical Institute of the Russian Academy of Sciences, Moscow, Russia}
 \email{zarhin\char`\@math.psu.edu}

\begin{document}
\begin{abstract}
A celebrated theorem of Bogomolov asserts that the $\ell$-adic Lie
algebra attached to the Galois action on the Tate module of an
abelian variety over a number field contains all homotheties. This
is not the case in characteristic $p$: a ``counterexample" is
provided by an ordinary elliptic curve defined over a finite
field. In this note we discuss (and explicitly construct) more
interesting examples of ``non-constant" absolutely simple abelian
varieties (without homotheties) over global fields in
characteristic $p$.
\end{abstract}

\maketitle

\section{Introduction}
Let $K$ be a field, $K_a$ its algebraic closure and
$\Gal(K)=\Aut(K_a/K)$ the absolute Galois group. If $X$ is an
abelian variety over $K$ then we write $\End_K(X)$ for the ring of
$K$-endomorphisms of $X$ and $\End_K^0(X)$ for the corresponding
$\Q$-algebra $\End_K(X)\otimes\Q$. We write $\End(X)$ for the ring
of $K_a$-endomorphisms of $X$ and $\End^0(X)$ for the
corresponding $\Q$-algebra $\End(X)\otimes\Q$. The notation $1_X$
stands for the identity automorphism of $X$. It is well-known
\cite{Mumford} that $\End^0(X)$ is a finite-dimensional semisimple
$\Q$-algebra and its center $\CC(X)$ is  a
product of number fields; in addition, either of those fields is
either totally real or a CM-field.

Let $E$ be a number field. Suppose we are given an embedding $$i:
E\hookrightarrow \End_K^0(X), \ i(1)=1_X.$$ Then $[E:\Q]$ divides
$2\dim(X)$ \cite[Ch. 2, Sect. 5, Prop. 2]{ShimuraA}; let us put
$$r(X,E)=\frac{2\dim(X)}{[E:\Q]}.$$

Let $\ell$ be a prime different from $\fchar(K)$. We write
$T_{\ell}(X)$ for the corresponding Tate $\Z_{\ell}$-module of $X$
and $V_{\ell}(X)$ for the corresponding $\Q_{\ell}$-vector space
$T_{\ell}(X)\otimes_{\Z_{\ell}}\Q_{\ell}$. It is well-known that
$T_{\ell}(X)$ is a free $\Z_{\ell}$-module of rank $2\dim(X)$ and
$V_{\ell}(X)$ is a $2\dim(X)$-dimensional $\Q_{\ell}$-vector
space. We write $\II$ for the identity automorphism of
$V_{\ell}(X)$. It is well-known that
$\Aut_{\Z_{\ell}}(T_{\ell}(X))\cong \GL(2\dim(X),\Z_{\ell})$ is a
compact $\ell$-adic Lie group with Lie algebra
$\End_{\Q_{\ell}}(V_{\ell}(X))$. Let
$$\det:\Aut_{\Q_{\ell}}(V_{\ell}(X))\to \Q_{\ell}^*$$
be the determinant map. As usual, we write $\SL(V_{\ell}(X))$ for
its kernel. It is well-known that $\SL(V_{\ell}(X))$ is a Lie
subgroup in $\Aut_{\Q_{\ell}}(V_{\ell}(X))$ and its Lie algebra
coincides with
$$\sL(V_{\ell}(X)):=\{u\in \End_{\Q_{\ell}}(V_{\ell}(X))\mid \tr(u)=0\}$$
where
$$\tr:\End_{\Q_{\ell}}(V_{\ell}(X)) \to \Q_{\ell}$$
is the trace map.

 On the other hand, $T_{\ell}(X)$ carries a natural
structure of $\End_K(X)\otimes\Z_{\ell}$-module and $V_{\ell}(X)$
carries a natural structure of
$\End_K^0(X)\otimes_{\Q}\Q_{\ell}$-module.

Let us put
$$E_{\ell}=E\otimes_{\Q}\Q_{\ell}\subset \End_K^0(X)\otimes_{\Q}\Q_{\ell}.$$
The embedding $i$ provides $V_{\ell}(X)$ with a natural structure
 of $E_{\ell}$-module: it is known \cite{Shimura,Ribet} that this
 module is free of rank $r(X,E)$.

 One may view $E_{\ell}^{*}$ as a commutative $\ell$-adic Lie (sub)group with (commutative) Lie
  algebra $E_{\ell}$. We have
  $$\Z_{\ell}^*\II\subset E_{\ell}^{*}\subset \Aut_{\Q_{\ell}}(V_{\ell}(X));$$
  clearly, $\Z_{\ell}^*\II$ is a compact $\ell$-adic Lie subgroup
  whose Lie algebra coincides with $\Q_{\ell}\II$.

  \begin{rem}
  \label{cartan}
Let $G \subset \Aut_{\Q_{\ell}}(V_{\ell}(X))$ be a compact
subgroup. Then the ($\ell$-adic variant of) Cartan's theorem
\cite[Part 2, Ch. 5, Sect. 9]{SerreL} tells us that $G$ is a Lie
subgroup. Clearly, the intersection $G\bigcap \Z_{\ell}^*\II$ is
infinite if and only if the Lie algebra $\Lie(G)$ of $G$ contains
$\Q_{\ell}\II$.
  \end{rem}

 Let us consider the centralizer $\End_{E_{\ell}}(V_{\ell}(X))$ of
 $E_{\ell}$ in $\End_{Q_{\ell}}(V_{\ell}(X))$ and its group of
 invertible elements $\Aut_{E_{\ell}}(V_{\ell}(X))$.
  One may view $\Aut_{E_{\ell}}(V_{\ell}(X))$ as  an
 $\ell$-adic Lie group with  Lie
 algebra $\End_{E_{\ell}}(V_{\ell}(X))$.

 Since
$V_{\ell}(X)$ is a free $E_{\ell}$-module of finite rank, there
are the natural $E_{\ell}$-determinant homomorphism of $\ell$-adic
Lie groups
$$\det{_{E_{\ell}}}:\Aut_{E_{\ell}}(V_{\ell}(X))\to E_{\ell}^*$$
and the $E_{\ell}$-trace map
$$\tr_{E_{\ell}}: \End_{E_{\ell}}(V_{\ell}(X))\to E_{\ell}.$$
Clearly, $\tr_{E_{\ell}}$ is the {\sl tangent map} of Lie algebras
attached to $\det{_{E_{\ell}}}$.

\begin{rem}
\label{finite} Let $G$ be a (closed) compact subgroup in
$\Aut_{E_{\ell}}(V_{\ell}(X))$. Then $G$ is an $\ell$-adic Lie
(sub)group and its Lie algebra $\Lie(G)$ is a $\Q_{\ell}$-Lie
subalgebra of $\End_{E_{\ell}}(V_{\ell}(X))$. In addition, if
$\Lie(G)$ is a semisimple Lie algebra then $\det{_{E_{\ell}}}(G)$
is a finite subgroup in $E_{\ell}^{*}$. Indeed, the semisimplicity
of $\Lie(G)$ implies that $\tr_{E_{\ell}}(\Lie(G))=\{0\}$ and
therefore $\det{_{E_{\ell}}}=1$ on an open subgroup of $G$. One
has only to recall that every open subgroup in a compact
$\ell$-adic Lie group has finite index.
\end{rem}

There is a natural continuous homomorphism ($\ell$-adic
representation) \cite{Serre}
$$\rho_{\ell,X}:\Gal(K) \to \Aut_{\Z_{\ell}}(T_{\ell}(X))\subset
\Aut_{\Q_{\ell}}(V_{\ell}(X));$$ its image $G_{\ell,X}$ is a
compact $\ell$-adic Lie subgroup of
$\Aut_{\Q_{\ell}}(V_{\ell}(X))$.  We write $\g_{\ell,X}$ for the
Lie algebra of $G_{\ell,X}$; one may view $\g_{\ell,X}$ as a Lie
$\Q_{\ell}$-subalgebra in $\End_{\Q_{\ell}}(V_{\ell}(X))$
\cite{Serre}.

 The
following assertion is proven in \cite{ZarhinT}.

\begin{thm}
\label{ZT}
 Suppose that $K$ is a global field of characteristic
$p>2$ and $X$ is an abelian variety of positive dimension over
$K$. Then:

\begin{itemize}
\item[(I)]
 $\g_{\ell,X}$ is a reductive
$\Q_{\ell}$ algebra, i.e. $\g_{\ell,X}\cong \g^{ss}\oplus \cc$
where $\g^{ss}$ is a semisimple $\Q_{\ell}$-Lie algera and $\cc$
is the center of $\g_{\ell,X}$.
\item[(II)]
 $\dim_{\Q_{\ell}}(\cc)=1$.
 \item[(III)]
 If $\CC(X)$ is  a product of
 totally real number fields then $\cc=\Q_{\ell}\cdot\II$.
 \end{itemize}
\end{thm}

When $K$ is a number field, a theorem of Bogomolov
\cite{Bogomolov, Bogomolov2} asserts that $\g_{\ell,X}$ always
contains homotheties $\Q_{\ell}\cdot\II$.

However, one may easily check that this is not the case if $K$ is
a global field of characteristic $p$. For example, if $X$ is an
ordinary elliptic curve that is defined over a finite field then
$\g_{\ell,X}$ is a one-dimensional $\Q_{\ell}$-Lie algebra that is
generated by the $\ell$-adic logarithm of the corresponding
Frobenius endomorphism, which is not a scalar. The aim of this
note is to prove the existence of an absolutely simple abelian
variety $X$ over a global field of characteristic $p$ such that
$\g_{\ell,X}$ does {\sl not} contain homotheties and $X$ is {\sl
not} isogenous over $K_a$ to an abelian variety over a finite
field. Recall \cite{Oort} that the latter condition means that $X$
is not an abelian variety of CM-type over $K_a$. Our main result
is described by the following two statements.

\begin{thm}
\label{odd}
 Suppose that $K$ is a global field of characteristic
$p>2$. Suppose that $X$ is an ordinary abelian variety of positive
dimension over $K$. Let $E \subset \End^0(X)$ be a subfield that
contains $1_X$. Assume that $r(X,E)$ is an odd integer.

Then $\g_{\ell,X}\bigcap \Q_{\ell}\cdot\II=\{0\}$, i.e.,
$\g_{\ell,X}$ does not contain homotheties except zero and
$G_{\ell,X}\bigcap \Z_{\ell}^*\II$ is finite.
\end{thm}

We prove Theorem \ref{odd} in Section \ref{OAV}.

\begin{thm}
\label{OPT} Let $Z$ be an ordinary elliptic curve over a finite
field $k$ of characteristic $p>2$ and $E=\End^0(Z)$ the
corresponding imaginary quadratic field.

Then for every odd $g>1$ there exist a global field $K$ of
characteristic $p$ and an ordinary $g$-dimensional abelian variety
$X$ over $K$ that enjoys the following properties:

\begin{itemize}
\item[(i)] All endomorphisms of $X$ are defined over $K$ and
$\End^0(X)=E$. In particular, $X$ is absolutely simple.
     \item[(ii)]
 $X$  is not
isogenous over $K_a$ to an abelian variety that is defined over a
finite field.
\item[(ii)] $\g_{\ell,X}\bigcap \Q_{\ell}\cdot\II=\{0\}$,
i.e., $\g_{ell,X}$ does not contain
homotheties except zero and $G_{\ell,X}\bigcap \Z_{\ell}^*\II$ is
finite.
\end{itemize}
\end{thm}

\begin{rem}
\begin{itemize}
\item[(i)]
  In light of Theorem 2(b) of \cite{TateInv},  the second assertion of Theorem
\ref{OPT} follows readily from the first one, because in this case
$$\dim_{\Q}(\End^0(X))=\dim_{\Q}(E)=2<2g=2\dim(X).$$
\item[(ii)]
In light of Theorem \ref{odd}, the third assertion of Theorem
\ref{OPT} follows readily from the first one, because in this case
$r(X,E)=g$ is odd.
\end{itemize}
\end{rem}

We prove Theorem \ref{OPT}(i) in Section \ref{MO}. In Section
\ref{super} we discuss an explicit example of an abelian variety
that satisfies the conditions and conclusions of Theorem
\ref{odd}.

I am grateful to B. Poonen, F. Voloch and M. Stoll for  a
stimulating question that was asked  during the special semester
``Rational and integral points on higher-dimensional varieties" at
the MSRI. My special thanks go to the MSRI and the organizers of
this program.  I am grateful to the referee, whose comments helped to
improve the exposition.

\section{Abelian varieties and imaginary quadratic fields}
\label{MO}

\begin{proof}[Proof of Theorem \ref{OPT}(i)]
Notice that all endomorphisms of $Z$ are defined over $k$.
 (This well-known result goes back to Deuring \cite{Deuring}; it follows easily from
 Main Theorem of \cite{TateInv}.)
Since
$Z$ is ordinary and $g-1$ is a multiple of $2=\dim_{\Q}(E)$,  a
theorem of Oort -van der Put \cite[Th. 1.1]{OP} implies the
existence of an ordinary $g$-dimensional abelian variety $Y$ over
$k((t))$ with all endomorphisms defined over $k((t))$  and
$\End^0(Y)=E$. Clearly, $Y$ and all its endomorphisms are defined
over a field $K$ that is finitely generated over $k$. Now, Mori's
specialization arguments \cite[Cor. 5.4]{Mori} allow us to assume
that $K$ has transcendence degree $1$, i.e., is global.
\end{proof}

\section{Ordinary abelian varieties}
\label{OAV}

\begin{lem}
\label{ordinary} Let $k$ be a finite field that consists of $q$
elements, $A$ an ordinary abelian variety over $k$ and $d$ a
positive odd integer. If $\{\alpha_1, \dots , \alpha_d\}$ are $d$
eigenvalues of the Frobenius endomorphism $Fr$ of $A$ then
$q^{-d}(\prod_{i=1}^d \alpha_i)^2$ is not a root of unity.
\end{lem}

\begin{proof}[Proof of Lemma \ref{ordinary}]
If $p=\fchar(k)$ then $q$ is a power of $p$. Let us choose a
 $p$-adic valuation map $\ord_p:\bar{\Q}^*\to \Q$ normalized by
 the condition $\ord_p(q)=1$. Since $A$ is ordinary,   the Honda-Tate theory
 \cite{T} tells us that
 $\ord_p(\alpha)=0$ or $1$ for every eigenvalue of the Frobenius endomorphism of
 $A$. This implies that
 $$\ord_p(q^{-d}(\prod_{i=1}^d \alpha_i)^2)=-d+2\sum_{i=1}^d
 \ord_p(\alpha_i)\in -d+2\Z$$
 is an {\sl odd} integer and therefore does not vanish. It follows
 that $q^{-d}(\prod_{i=1}^d \alpha_i)^2$ is not a root of unity.
\end{proof}

\begin{proof}[Proof of Theorem \ref{odd}]
Replacing (if necessary) $K$ by its finite separable algebraic
extension, we may and will assume that all endomorphisms of $X$
are defined over $K$; in particular, $E\subset
\End_K^0(X)=\End^0(X)$. Let us assume that $\g_{\ell,X}\bigcap
\Q_{\ell}\cdot\II\ne \{0\}$. This means that $\g_{\ell,X}$
contains $\Q_{\ell}\cdot\II$ and therefore
$\cc=\Q_{\ell}\cdot\II$.

Let us put $G^0=G_{\ell,X}\bigcap SL(V_{\ell}(X)$. Clearly, $G^0$
is a closed (compact) Lie subgroup of $G_{\ell,X}$ and $\Lie(G^0)$
has codimension $1$ in $\Lie(G_{\ell,X})=\g^{ss}\oplus
\Q_{\ell}\cdot\II$. The semisimplicity of $\g^{ss}$ implies
that $\Lie(G^0)=\g^{ss}$.

Let us put $$S=G_{\ell,X}\bigcap (1+\ell^2\Z_{\ell})\II\subset
\Z_{\ell}^*\II.$$ Clearly, $S$ is compact. Since
$\g_{\ell,X}=\Lie(G_{\ell,X})$ contains
 $\Q_{\ell}\cdot\II$, the group $G_{\ell,X}$ contains an open
 subgroup of $\Z_{\ell}^*\II$. It follows that $S$ is an open
 subgroup of  finite index in $\Z_{\ell}^*\II$. Since
 $1+\ell^2\Z_{\ell}$ does not contain nontrivial roots of unity,
  $S$ does not contain elements of finite order (except $\II$) and
  therefore $G^0\bigcap S=\{\II\}$. Recall that both $G^0$ and $S$
  are subgroups of $G_{\ell,X}$. Let us consider the homomorphism
  of compact $\ell$-adic Lie groups
  $$\pi: G^0\times S \to G_{\ell,X}, \ (u,c)\mapsto uc=cu.$$
  Clearly, $\pi$ is injective and the corresponding tangent map of Lie algebras is
  an isomorphism. It follows that $G^1:=\pi(G^0\times S)$ is an
  open
  compact subgroup in $G_{\ell,X}$ and $\pi$ induces an
  isomorphism of $\ell$-adic Lie groups $G^0\times S$ and $G^1$.

  \begin{lem}
  \label{G1}
There exists a positive integer $m$ such that
$${\det}_{E_{\ell}}(g)^m\in \Q_{\ell}^*\II \ \forall g \in
G_{\ell,X}.$$
\end{lem}

\begin{proof}[Proof of Lemma \ref{G1}]
Since $\Lie(G^0)$ is semisimple, it follows from Remark
\ref{finite} that ${\det}_{E_{\ell}}(G^0)$ is a finite group. If
$m_0$ is its order then ${\det}_{E_{\ell}}(g_0)^{m_0}=1$ for all
$g_0\in G^{0}$. Notice that
$${\det}_{E_{\ell}}(c)=c^{r(X,E)} \ \forall c \in \Z_{\ell}^*\II,$$ because
$\Z_{\ell}^*\II\subset E_{\ell}^*$. It follows that
${\det}_{E_{\ell}}(g)^{m_0}\in \Q_{\ell}^*\II \ \forall g \in
G^{1}$. In order to finish the proof, one has only to recall that
$G^1$ is a subgroup of finite index in $G_{\ell,X}$ and put
$m:=m_0\cdot [G_{\ell,X}:G^1]$.
\end{proof}

There exists a place $v$ of $K$ such that the abelian variety $X$
has ordinary good reduction. (In fact, this condition is fulfilled
for all but finitely many places of $K$.) Let $k(v)$ be the
residue field at $v$, let $q(v)$ be the cardinality of $k(v)$ and
$X(v)$ the reduction of $X$ at $v$, which is an ordinary abelian
variety over $k(v)$ whose dimension coincides with $\dim(X)$.  Let
$\P_v(t)\in \Z[t]$ be the (degree $2\dim(X)$) characteristic
polynomial of the Frobenius endomorphism $\Fr$ of $X(v)$. One may
view the roots of $\P_v$ as eigenvalues of the Frobenius endomorphism
with respect to its natural action on $V_{\ell}(X(v))$.

Let us choose a place $\bar{v}$ of $K_a$ that lies above $v$. Such a choice gives rise
to natural isomorphisms \cite{SerreTate,Serre}
$$T_{\ell}(X)\cong T_{\ell}(X(v)), \ V_{\ell}(X)\cong
V_{\ell}(X(v))$$ in such a way that $\Fr\in
\Aut_{\Z_{\ell}}(T_{\ell}(X(v))$ corresponds to a certain element
of $G_{\ell,X}$: this element is called the {\sl Frobenius
element} attached to $\bar{v}$ and denoted by $F_{\bar{v}}$. It is
known \cite[Chap. 7, proof of Prop. 7.23]{Shimura} (see also
\cite[p. 167]{ZarhinInv}) that
$$b_v:={\det}_{E_{\ell}}(F_{\bar{v}})\in E^*\subset E_{\ell}^*$$
and $b_v$ is a product of $r(X,E)$ eigenvalues of $\Fr$.

In other words, let $L$ be the splitting field of $\P_v(t)$ over
$E$: it is a finite Galois extension of $E$. Then there exist
roots $\alpha_1, \dots , \alpha_{r(X,E)}$ of $\P_v(t)$ such that
their product coincides with $b_v$. On the other hand, it follows
from a famous theorem of A. Weil (the Riemann hypothesis)
\cite[Sect. 21]{Mumford} that if we fix a field embedding
$L\subset \C$ then
$$\mid b_v^2\mid_{\infty}=q(v)^{r(X,E)}$$
where $\mid \ \mid_{\infty}$ is the standard (archimedean)
absolute value on the field of complex numbers.
 On the other hand, by Lemma \ref{G1}, there exists a
positive integer $m$ such that $b_v^{m}\in \Q_{\ell}$. Since the
intersection of $E=E\otimes 1$ and $\Q_{\ell}=1\otimes \Q_{\ell}$
in $E_{\ell}=E\otimes_{\Q}\Q_{\ell}$ coincides with $\Q$, we
conclude that $b_v^{m}$ is a rational number. This implies that
$b_v^{2m}$ is a positive rational number and, by Weil's theorem,
coincides with $q(v)^{m r(X,E)}$. This implies that
$$1=\left(q(v)^{-r(X,E)}\cdot b_v^2\right)^m.$$
However, by Lemma \ref{ordinary}, $q(v)^{-r(X,E)}b_v^2$ is {\sl
not} a root of unity. (Here we use the {\sl oddity} of $r(X,E)$.)
We get a contradiction, which proves the Theorem.
\end{proof}

\section{Superelliptic jacobians}
\label{super}

\begin{prop}
\label{cmord}
 Let $K$ be a number field with the ring
of integers $\O_K$. Let $Y$ be an abelian variety of positive
dimension over $K$, let $L$ be a CM-field of degree $2\dim(X)$ and
$i: L\hookrightarrow \End^0(Y)$ an embedding that sends $1$ to
$1_Y$. Let $p$ be a prime that splits completely in $L$, i.e.
$L\otimes_{\Q}\Q_p$ splits into a product of $[L:\Q]$ copies of
$\Q_p$. Let $\p$ be maximal ideal in $\O_K$ with residual
characteristic $p$.

If $Y$ has good reduction at $\p$ then this reduction is ordinary.
\end{prop}

\begin{proof}
Let $\bar{\Q}_p$ be an algebraic closure of $\Q_p$. Let $L_{\p}$
be the $\p$-adic completion of $L$. By assumption, $L_{\p}=\Q_p$
and therefore the set $H_{\p}$ of $\Q_p$-linear field embeddings
$L_{\p}\hookrightarrow \bar{\Q}_p$ is a singleton that consists of
the inclusion map $\Q_p\subset \bar{\Q}_p$; in particular,
$\#(H_{\p})=1$. Now the assertion follows readily from Lemma 5 in
Sect. 4 of \cite{T}.
\end{proof}

\begin{lem}
\label{t0}
Let us consider the curve $C_0:y^3=x^9-x$ and its
jacobian $J(C_0)$ over $\Q$.

Then:
\begin{itemize}
 \item [(i)]If $p$ is a prime such that $p-1$ is divisible by $24$
then $J(C_0)$ has ordinary good reduction at $p$. \item[(ii)]
$J(C_0)$ is a (non-simple) abelian variety of CM-type over
$\bar{\Q}$.
\end{itemize}
\end{lem}

\begin{proof}
Clearly, both $C_0$ and $J(C_0)$ have good reduction at $p$,
because $x^9-x=x(x^8-1)$ has $9$ distinct roots in $\F_p$ and
therefore has no multiple roots in characteristic $p$. In order to
check that $J(C_0)$ has ordinary reduction,
 pick a number field $F$ such that $F$ contains $\Q(\zeta_{24})$,
 all endomorphisms of $J(C_0)$  are defined over $F$
 and all homomorphisms between $J(C_0)$ and the
 elliptic curve $y^2=x^3-x$ are defined over $F$. Let us
 consider both $C_0$ and $J(C_0)$ over $F$, and let $\p$ be a place of $F$ that lies above
$p$. For our purposes, it suffices to check that $J(C_0)$ has
ordinary reduction at $\p$.

Pick  a primitive cubic root of unity $\zeta_3 \in F$. Then the
map
$$(x,y)\mapsto (x,\zeta_3 y)$$
induces an automorphism $\delta_3:C_0 \to C_0$, which, in turn,
induces by Albanese functoriality an automorphism $J(C_0)\to
J(C_0)$, which we still denote by $\delta_3$. It is known \cite[p.
149]{Poonen} that $\delta_3^2+\delta_3+1=0$ in $\End(J(C_0))$,
which leads to the embedding
$$\Z[\zeta_3]\hookrightarrow \End(J(C_0)), \
\zeta_3\mapsto\delta_3, 1\mapsto 1_{J(C_0)}.$$ Extending it by
$\Q$-linearity, we get an embedding
$$\Q(\sqrt{-3})=\Q(\zeta_3)\hookrightarrow \End^0(J(C_0)), \
\zeta_3\mapsto\delta_3, 1\mapsto 1_{J(C_0)}.$$ On the other hand,
pick a primitive $8$th root of unity  $\zeta_8\in F$. Then the map
$$(x,y)\mapsto (\zeta_8^{-1}x,\zeta_8^{-3} y)$$
induces an automorphism $\delta_8:C_0 \to C_0$, which commutes
with $\delta_3$. Again $\delta_8$  induces by Albanese
functoriality an automorphism  of $J(C_0)$, which we still denote
by $\delta_8$; clearly $\delta_8$ and $\delta_3$ do commute in
$\End(J(C_0)$. In order to understand $\delta_8$ better, let us
divide both sides of the equation for $C_0$ by $x^9=(x^3)^3$: we
get $(y/x^3)^3=1-(1/x)^8$. It follows that $C$ is $F$-birationally
isomorphic to the curve $$C^{\prime}: w^8=-u^3+1; \ w=1/x,
u=y/x^3$$ and $\delta_8$ is induced by
$$(u,w)\mapsto (u,\zeta_8 w).$$
This implies that the jacobian $J(C^{\prime})$ of $C^{\prime}$ and
$J(C)$ are isomorphic over $F$. Let us put $f(w)=-w^3+1$. Then in
notations of \cite{ZarhinP}, $C^{\prime}=C_{f,8}$ and the
structure of its jacobian $J(C^{\prime})=J(C_{f,8})$ is described
as follows \cite[Sect. 5, Cor. 5.12, Rem. 5.14, Th. 5.17
]{ZarhinP}. First, $J(C_{f,8})$ contains a $\delta_8$-invariant
abelian fourfold
$$J^{(f,8)}=(\delta_8^3+\delta_8^2+\delta_8+1)(J(C_{f,8}))\subset J(C_{f,8})$$
provided with an embedding
$$\Z[\zeta_8]\hookrightarrow \End(J^{(f,8)}), \ \zeta_8\mapsto
\delta_8, 1 \mapsto 1_{J(C_{f,8})}.$$ Clearly, $J^{(f,8)}$ is
$\delta_3$-invariant. This gives rise to an embedding
$$\Q(\zeta_{24})=\Q(\zeta_3)\otimes\Q(\zeta_8)\hookrightarrow
\End^0(J^{(f,8)}), \ \zeta_3\mapsto \delta_3, \zeta_8\mapsto
\delta_8$$ and $1$ goes to the identity map. This implies that
$J^{(f,8)}$ is an abelian fourfold of CM-type. Since $p$ splits in
$\Q(\zeta_{24})$,  it follows from Proposition \ref{cmord} that
$J^{(f,8)}$ has ordinary reduction at all places of $F$ over $p$.
Second, $J(C_{f,8})$ is isogenous (over $F$) to a product of
$J^{(f,8)}$, two copies of the elliptic curve $y^2=x^3-x$ and the
elliptic curve $w^2=-v^3+1$. Since $24$ divides $p-1$, the prime
$p$ splits in the imaginary quadratic fields $\Q(\sqrt{-1})$ and
$\Q(\sqrt{-3})$. Therefore the CM-elliptic curves $y^2=x^3-x$ with
multiplication by $\Q(\sqrt{-1})$ and $w^2=-v^3+1$ with
multiplication by $\Q(\sqrt{-3})$  have ordinary reduction at $p$.
It follows that $J(C_{f,8})$ has ordinary reduction at $\p$.
\end{proof}

\begin{ex}
Fix a prime $p$ with $p-1$ divisible by $24$, let $K=\F_p(t)$ and
let $X$ be the $7$-dimensional jacobian of the $K$-curve
$C:y^3=x^9-x-t$. Since $p$ divides neither $9$ nor $8$, $x^9-x\in
\F_p[x]$ is a {\sl Morse polynomial} \cite[p. 39]{SerreG}, i.e.,
its derivative $9 x^8-1$ has $8$ distinct roots $\beta_1, \dots,
\beta_8$ and all eight critical values
$\beta_i^9-\beta_i=-\frac{8}{9}\beta_i$ are distinct. It follows
that the Galois group of $x^9-x-t$ over $\F_p(t)$ is the full
symmetric group $\ST_9$ \cite[p. 41]{SerreG}. On the other hand,
if $\zeta \in \F_p$ is a primitive cubic root of unity then
$$(x,y)\mapsto (x,\zeta y)$$
gives rise to a non-trivial automorphism of $C$ (of period $3$), which, in turn,
allows us to define the embedding
$$\Q(\zeta_3)=\Q(\sqrt{-3})\hookrightarrow \End^0(J(C)), \ 1 \mapsto 1_{J(C)}.$$
By Theorem 0.1 of \cite{ZarhinSb}, $E=\Q(\sqrt{-3})$ coincides
with its own centralizer in $\End^0(J(C))$ and therefore contains
$\CC(J(C))$. This means that $\CC(J(C))=E$ or $\Q$. On the other
hand,  the reduction of $J(C)$ at $t=0$ is the jacobian of the
$\F_p$-curve $y^3=x^9-x$, which is ordinary, by Lemma \ref{t0}.
Applying Theorem \ref{odd}, we obtain that $\g_{\ell,J(C)}$ does
not contain non-zero homotheties. On the other hand, if
$\CC(J(C))=\Q$ then, by Theorem \ref{ZT}(III), $\g_{\ell,J(C)}$
does  contain all the homotheties. This contradiction proves that
$\CC(J(C))=E$ and therefore the centralizer of $E$ coincides with
the whole $\End^0(J(C))$. This implies that $\End^0(J(C))=E$ and
therefore $J(C)$ is absolutely simple and  is not of CM-type. It
follows that $J(C)$ is not isogenous to an abelian variety that is
defined over a finite field.
\end{ex}


\begin{thebibliography}{99}

\bibitem{Bogomolov} F. A. Bogomolov, {\em Sur l{'}alg\'ebricit\'e des repr\'esentations}
 $\ell$-{\em adiques}. C.R. Acad. Sci. Paris S\'er. A-B (1980), no. 15, A701--A703.
 \bibitem{Bogomolov2} F. A. Bogomolov, {\em Points of finite order on an abelian
 variety}. Izv. Akad. Nauk SSSR Ser. Mat. {\bf 44} (1980), 782--804; Math.
USSR Izv. {\bf 17} (1981), 55--72.

\bibitem{Deuring} M. Deuring, {\em Die Typen der Multiplikatorenringe elliptischer
Funktionenk\"orper}. Abh. Math. Sem. Hamburg {\bf 14} (1941), 197--272.

\bibitem{Mori} Sh. Mori, {\em On Tate conjecture concerning endomorphisms of abelian varieties}.
 Proceedings of the International Symposium on Algebraic Geometry
(Kyoto Univ., Kyoto, 1977), pp. 219--230, Kinokuniya Book Store,
Tokyo, 1978.

\bibitem{Mumford} D. Mumford,  Abelian varieties, 2nd edn, Oxford University Press, 1974.

\bibitem{Oort} F. Oort, {\em The isogeny class of a CM-abelian
variety is defined over a finite extension of the prime field}. J.
Pure Applied Algebra {\bf 3} (1973), 399--408.


\bibitem{OP} F. Oort, M.-van der Put, {\em A construction of an abelian variety with a given
endomorphism algebra}. Compositio Math. {\bf 67} (1988), 103--120.

\bibitem{Poonen} B. Poonen and E. Schaefer, {\em Explicit descent for Jacobians of cyclic covers of
 the projective line}. J. reine angew. Math. {\bf 488} (1997), 141--188.

\bibitem{Ribet} K. Ribet, {\em Galois action on division points of Abelian varieties with real
multiplications}. Amer. J. Math. {\bf 98} (1976), 751--804.

\bibitem{Serre} J.-P. Serre, Abelian $\ell$-adic representations
and elliptic curves, 2nd edition, Addison Wesley, 1989.

\bibitem{SerreG} J.-P. Serre, Topics in Galois Theory,
Jones and Bartlett Publishers, Boston-London, 1992.

\bibitem{SerreTate} J-P.\ Serre, J.\ Tate,
{\em  Good reduction of abelian varieties}. Ann.\ of Math.\ (2)
{\bf 88}  (1968), 492--517.

\bibitem{SerreL} J.-P. Serre, Lie algebras and Lie groups, 2nd
edition, Springer Lecture Notes in Math. {\bf 1500} (1992).

\bibitem{Shimura} G.  Shimura, Introduction to the arithmetic
theory of automorphic functions. Publ. Math. Soc. Japan {\bf 11},
Princeton University Press, 1971.

\bibitem{ShimuraA} G. Shimura, Abelian varieties with complex
multiplication and modular functions. Princeton University Press,
Princeton, 1997.

\bibitem{T} J. Tate, {\em Classes d'isog\'enie des vari\'etes ab\'eliennes sur un corps fini
(d'apr\'es Honda)}.  S\'eminaire Bourbaki {\bf 352} (1968).
Springer Lecture Notes in Math. {\bf 179} (1971), 95--110.

\bibitem{TateInv} J. Tate, {\em   Endomorphisms of abelian varieties over
 finite fields}. Inv. Math. {\bf 2} (1966), 134--144.

\bibitem{ZarhinT} Yu. G. Zarhin, {\em Torsion of abelian varieties
in finite characteristic}.  Mat. Zametki {\bf 22} (1977), 3--11;
Math. Notes {\bf 22} (1978),  493--498.

\bibitem{ZarhinInv} Yu. G. Zarhin, {\em Abelian varieties},
$\ell$-{\em adic representations and Lie algebras. Rank
independence on} $\ell$. Inv. Math. {\bf 55} (1979), 165 - 176.

\bibitem{ZarhinSb} Yu. G. Zarhin, Endomorphism rings of certain jacobians in finite characteristic.
 Matem. Sbornik  {\bf 193} (2002), issue 8, 39--48;
 Sbornik Math., 2002, {\bf 193} (8), 1139-1149.

\bibitem{ZarhinP} Yu. G. Zarhin, {\em Superelliptic jacobians}.
arXiv:math.AG/0601072; to appear in Proceedings of the Pisa
research programme on Diophantine Geometry (Spring 2005).

\end{thebibliography}
\end{document}